\documentclass[12pt]{article}
\usepackage{amssymb,amsfonts,amsmath,geometry,esint}

\catcode`\@=11
\@addtoreset{equation}{section}

\catcode`\@=12
\usepackage{color}

\def\R{\mathbb{R}}

\def\f{\varphi}

\def\irn{\int\limits_{\R^n}}
\def\iirn{\int\limits_0^\infty\!\int\limits_{\R^n}}

\def\eps{\varepsilon}
\def\div{{\rm div}}

\def\lh{{\lambda_h}}

\def\Ds{\left(-\Delta\right)^{\!s}\!} %Stanfard Dirichlet-Fourier
\def\Dshalf{\left(-\Delta\right)^{\!\frac {s}{2}}\!}  %Stanfard Dirichlet-Fourier s/2

\def\S{\mathbb S}

\def\sstar{{2^*_s}}

\def\bu{{z}}
\def\budot{\dot{\bu}}
\def\buddot{\ddot{\bu}}

\def\bq{b}

\def\proof{\noindent{\textbf{Proof. }}}
\def\QED{\hfill {$\square$}\goodbreak \medskip}

\newtheorem{Theorem}{Theorem}[section]
\newtheorem{Lemma}[Theorem]{Lemma}
\newtheorem{Proposition}[Theorem]{Proposition}
\newtheorem{Corollary}[Theorem]{Corollary}
\newtheorem{Remark}[Theorem]{Remark}

\usepackage{mathtools}

\begin{document}

%\title {The equations related to the fractional Hardy--Sobolev inequality: \\ qualitative properties of solutions}

\title {Complete classification and nondegeneracy\\ of minimizers for the fractional \\ Hardy-Sobolev inequality, and applications}

\author{Roberta Musina\footnote{Dipartimento di Matematica ed Informatica, Universit\`a di Udine, via delle Scienze, 206 -- 33100 Udine, Italy. Email: {roberta.musina@uniud.it}. Partially supported by PRID-DMIF Project VAPROGE, Universit\`a di Udine.}~ 
and \setcounter{footnote}{6}
Alexander I. Nazarov\footnote{
St.Petersburg Dept of Steklov Institute, Fontanka 27, St.Petersburg, 191023, Russia, 
and St.Petersburg State University, 
Universitetskii pr. 28, St.Petersburg, 198504, Russia. E-mail: al.il.nazarov@gmail.com. Partially supported by RFBR grant 20-01-00630.
}
}

\date{}

\maketitle

\begin{abstract}
We study linear and non-linear equations related to the fractional Hardy--Sobolev inequality. We prove nondegeneracy of ground state solutions to the basic equation and investigate existence and qualitative properties, including symmetry of solutions to some perturbed equations.
\end{abstract}

%\noindent
%{\small {\bf Abstract.} }


\medskip\noindent
{\small {\bf 2010 Mathematics Subject Classfication:}}  35R11; 35B20; 35P99.
%\end{abstract}

\normalsize

\bigskip

\section{Introduction}
Our starting point is the nonlocal problem
\begin{equation*}\tag{$\mathcal P_0$}
\label{eq:CKN}
\begin{cases}
\Ds u=|x|^{-bq}u^{q-1}~,\quad  u\in \mathcal D^s(\R^n)\\
u>0~\!.
\end{cases}
\end{equation*}
Here $n\ge 2$, $\Ds$ is the fractional Laplacian of order $s\in(0,1)$, the exponents 
$q, b$ satisfy
\begin{equation}
2<q<\sstar:=\frac {2n}{n-2s}~,\qquad 
\frac {n}{q}-b=\frac {n}{2}-s
\label{eq:bq}
\end{equation}
and $\mathcal D^s(\R^n)$ is the natural Sobolev-type function space.
Problem (\ref{eq:CKN}) is  related to the fractional Hardy--Sobolev inequality
\begin{equation}
\label{eq:HS_ineq}
S_q\cdot\||x|^{-b} u\|_q^2\le\|\Dshalf u\|_2^2,\qquad u\in\mathcal D^s(\R^n),
\end{equation}
that plainly follows via H\"older interpolation between the Hardy and Sobolev inequalities. 

The best constant $S_q$ in (\ref{eq:HS_ineq}) is attained by a nonnegative radially symmetric function 
$$
\bu_1\in \mathcal D^s(\R^n)
$$
(see  \cite{MN_SB}) which is a weak solution to  $\Ds u=|x|^{-bq}u^{q-1}$.  
Since %we have
$\Ds \bu_1\ge 0$  in the sense of distributions, then the strong maximum principle 
(see \cite[Section 2]{Sil} and  \cite[Corollary 4.2]{MaxPri}), ensures that 
$\bu_1$  is lower semicontinuous and positive on $\R^n$. Hence, $\bu_1$ solves (\ref{eq:CKN}).
Further, by adapting  the moving plane argument in \cite{CLO} or \cite{DMPS} one can prove that $\bu_1$
 is radially symmetric about the origin and radially decreasing.

 We agree that the minimizer $\bu_1$ is fixed, form now on.
 
By direct computations one can check that for any $t>0$, the radial function 
$$
\bu_t(x)= t^{\frac {2s-n}{2}}\bu_1\big(\frac{x}{t}\big)
$$
achieves $S_q$ and solves (\ref{eq:CKN}). 
However, we emphasise the fact that, differently from the critical case $q=\sstar$  and from the local case $s=1$
(see \cite{CoTa,Li}, respectively),
the minimizers for $S_q$ are not explicitly known, nor classified. 

We are in position to state our first main result.

\begin{Theorem}[Regularity, decay estimates and uniqueness]
\label{T:CKN}
~
\begin{itemize}
\item[$i)$] $\bu_1\in {\cal C}^\infty(\R^n\setminus \{0\})\cap L^\infty(\R^n)\cap C^\alpha(\R^n)$ for any $\alpha\in[0, 2s-bq)$; moreover, there exist positive constants $C_1, C_2$ such that 
$$\frac{C_1}{1+|x|^{n-2s}}\le \bu_1(x)\le \frac{C_2}{1+|x|^{n-2s}}\qquad \text{for any $x\in\R^n$;}
$$
\item[$ii)$] if $u\in \mathcal D^s(\R^n)$ is a solution to {\em (\ref{eq:CKN})} then $u=\bu_t$ for some $t>0$;
\item[$iii)$] the function $t\mapsto \bu_t$  is a regular curve in $\mathcal D^s(\R^n)$ of class {${\cal C}^2$}.
\end{itemize}
\end{Theorem}

The proof of Theorem \ref{T:CKN} is based on some preliminary results  on  eigenvalue problems of the form
\begin{equation}
\label{eq:general_eigenvalue}
\Ds \f=\mu V(x) \f~,\quad \f\in \mathcal D^s(\R^n),
\end{equation}
where $V>0$ is a given measurable weight satisfying 
suitable integrability assumptions. Our results on (\ref{eq:general_eigenvalue}), see  Section \ref{S:eigenvalue},  might have an independent interest.
The proof
of Theorem \ref{T:CKN} is carried out in Section \ref{S:T11}.

Our next focus is the problem
\begin{equation*}\tag{$\mathcal L_t$}
\label{eq:lin}
\Ds v=(q-1)|x|^{-bq}\bu_t^{q-2}v~\!,\qquad v\in \mathcal D^s(\R^n),
\end{equation*}
which is obtained by linearizing (\ref{eq:CKN}) at $\bu_t$. Let us denote by 
a "dot" the differentiation with respect to $t$. Thanks to part $iii)$ in Theorem (\ref{T:CKN}),
%we can differentiate the identity
%$\Ds \bu_t=|x|^{-bq}\bu_t^{q-1}$
%with respect to $t$, to obtain 
it is easily seen
that $\budot_t$ is a weak solution to (\ref{eq:lin}). In Section \ref{S:T12} we prove 
 the next uniqueness result.

\begin{Theorem}[Nondegeneracy]
\label{T:nondegenerate}
If a function $v\in \mathcal D^s(\R^n)$ solves (\ref{eq:lin}), then $v$ is proportional to $\budot_t$.
\end{Theorem}

Nondegeneracy in the limiting case $q=\sstar$ has been proved in \cite{DDS}, by 
taking advantage of the explicit knowledge of the minimizer $\bu_t$.

As a  first consequence of Theorem \ref{T:nondegenerate}
 we  obtain a symmetry result for ground state (i.e. least energy) solutions to the nonlocal problem
\begin{equation*}\tag{$\mathcal P_\lambda$}
\label{eq:l_CKN}
\begin{cases}
\Ds u+\lambda|x|^{-2s}u=|x|^{-bq}u^{q-1}\qquad u\in \mathcal D^s(\R^n)\\
u>0~\!.
\end{cases}
\end{equation*}
If $\lambda\le 0$, then 
the moving plane method can be applied to show that any weak solution to (\ref{eq:l_CKN}) is radially symmetric about the origin. In particular,
letting $H_s$ to be the fractional Hardy constant (see \cite{He} for its explicit value), we have that 
any minimizer for the best constant
\begin{equation}
\label{eq:lambda}
S^\lambda_q=\inf_{u\in\mathcal D^s(\R^n)\atop u\neq 0} \frac {\|\Dshalf u\|_2^2+\lambda\||x|^{-s} u\|_2^2}{\||x|^{-b} u\|_q^2}
\end{equation}
is radial, provided that $-H_s<\lambda\le 0$ (existence has been proved in \cite{MN_SB}). On the other hand,
symmetry breaking occurs: 
if $\lambda>0$ is large, then  no extremal for $S^{\lambda}_q$ is radially symmetric
(see \cite[Theorem 1.1]{MN_SB}). 

In the next theorem, which is proved in Section \ref{S:positivity}, we show that symmetry 
persists also for small positive values of $\lambda$.

\begin{Theorem}[Symmetry preserving]
\label{T:symmetry}
There exists $\lambda^{\!R}_{s}>0$ such that for every 
$\lambda\in(-H_s,\lambda^{\!R}_{s})$, any minimizer for $S^{\lambda}_q$ is radially symmetric about the origin. 
\end{Theorem}

As a further consequence of Theorem \ref{T:nondegenerate}, in Section \ref{S:existence} we use a Lyapunov-Schmidt argument
inspired by  \cite[Sections 3 and 4]{FeSc} to obtain
sufficient conditions on  a prescribed weight $k(x)$ on $\R^n$ which guarantee the existence of solutions to the  perturbative model problem 
\begin{equation*}\tag{$\mathcal P^\eps_k$}
\label{eq:eps_CKN}
\begin{cases}
\Ds u=(1+\eps k(x))|x|^{-bq}u^{q-1}\qquad u\in \mathcal D^s(\R^n)\\
u>0~\!.
\end{cases}
\end{equation*}
For instance, we obtain the following extension of \cite[Theorem 1.3]{FeSc}. % to the fractional setting.

\begin{Theorem}
\label{T:pert0}
Let $k\in L^\infty(\R^n)$.
If $\lim\limits_{x\to 0}k(x)= \lim\limits_{|x|\to \infty}k(x)$, then problem {\em (\ref{eq:eps_CKN})} has at least a solution
for any $\eps$ close enough to $0$.
\end{Theorem}

%
%Our paper has the following structure. 
%
%Theorems \ref{T:CKN} and \ref{T:nondegenerate} are proved in Sections \ref{S:T11} and \ref{S:T12} respectively.
%The proof of Theorem \ref{T:pert0} is obtained by adapting the arguments in \cite{FeSc}, see Section \ref{S:existence}. 

%\medskip

{\small 
\paragraph{Notation.} The fractional Laplacian $\Ds$ in $\R^n$, $n\ge 2$, is formally defined by 
$$
{\mathcal F}\big[\Ds u\big] = |\xi|^{2s}{\mathcal F}[u]~\!,
$$
where ${\mathcal F}=\displaystyle{{\mathcal F}[u](\xi)= (2\pi)^{-\frac{n}{2}}
\int_{\R^n} e^{-i~\!\!\xi\cdot x}u(x)~\!dx}$ is the Fourier transform.
Thanks to the Sobolev inequality, the space
$$
\mathcal D^s(\R^n)=\big\{u\in L^\sstar(\R^n)~|~\Dshalf u \in L^2(\R^n)~\!\big\}
$$
naturally inherits a Hilbertian structure from the relations
$$
(u,v)_{\mathcal D^s}=(\Ds u,v)=\irn\Dshalf u\Dshalf v~\!dx~\!, \qquad
\|u\|_{\mathcal D^s}^2=(u,u)_{\mathcal D^s}. %=\irn|\xi|^{2s}\mathcal F[u]\overline{\mathcal F[v]}~\!d\xi~\!.
$$
From now on, we will always use the shorter notation $\mathcal D^s$ instead of $\mathcal D^s(\R^n)$, and we let $(\mathcal D^s)'$ be its dual space.
By elementary arguments, any $w\in (\mathcal D^s)'$ can be identified with the distribution $\Ds v$, where $v\in\mathcal D^s$ is uniquely determined by $w$.

Denote by $\|\cdot\|_p$ the norm in $L^p(\R^n)$.

For $0<\alpha<1$, ${\cal C}^\alpha$ stands for standard H\"older space. For $1<\alpha<2$, we denote by ${\cal C}^\alpha$ the space of continuously differentiable functions with $\nabla u\in {\cal C}^{\alpha-1}$.
}

%\red{{Constants depending only on $n, q, \|k\|_\infty$ %and possibly on the choice of $\bu_1$ 
%are denoted by $c$,
%apart from few constants that play a special role.}}

\section{Preliminaries on  eigenvalue problems}
\label{S:eigenvalue}

In this section we study the linear problem (\ref{eq:general_eigenvalue}) under the assumption
$V>0$.
We use the following regularity results within the classical theory for Riesz potentials.

\begin{Proposition}\label{Riesz}
Let  $\alpha\in(0,n)$ be given.

 $i)$ Let $f\in L^p(\R^n)$ with $p<\frac n{\alpha}$, then $f*|x|^{\alpha-n}\in L^{\frac {np}{n-\alpha p}}(\R^n)$ (the Hardy--Littlewood--Sobolev theorem, \cite[Ch. V, Theorem 1]{St});

$ii)$ Let $f\in L^p(\R^n)$ with $p>\frac n{\alpha}$; if  $\alpha>1$ assume in addition that $p<\frac n{\alpha-1}$. Then $f*|x|^{\alpha-n}\in {\cal C}^{\alpha-\frac np}(\R^n)$ (\cite[Ch. V, Theorem 5]{St} and \cite[Ch. V, 6.7a)]{St}).
\end{Proposition}

We say that a nontrivial function $\f\in\mathcal D^s$ is 
an eigenfunction for (\ref{eq:general_eigenvalue}) if it is a weak solution to (\ref{eq:general_eigenvalue}), namely
$$
(\Ds\f,\psi)=\mu\irn V(x)\f\psi~\!dx\quad\text{for any $\psi\in{\cal C}^\infty_0(\R^n)$.}
$$

\begin{Lemma}
\label{L:general_spectrum} 1. Let $V\in L^{\frac n{2s}}(\R^n)$. Then the spectrum of (\ref{eq:general_eigenvalue}) is discrete. We denote by $\mu_j$ a non-decreasing 
unbounded
sequence of  eigenvalues counting with multiplicities. The corresponding eigenfunctions $\f_j$ form a complete orthogonal system in $\mathcal D^s$. Moreover, 
\begin{equation}
\label{eq:eigen}
\mu_j=\min_{\scriptstyle
\f\in\mathcal D^s,\  \f\neq 0 \atop 
\scriptstyle(\f,\f_i)_{\mathcal D^s}=0,\ \forall i<j}\ 
\frac{\displaystyle(\Ds \f,\f)}{\displaystyle\irn V(x)|\f(x)|^2~\!dx}~\!.
\end{equation}
The first eigenvalue $\mu_1>0$ is simple, and it is the only eigenvalue admitting a positive eigenfunction.
 
2. If in addition  $V\in L^{\frac n{2s-\eps}}(\Omega)$ for some $\eps\in(0,2s)$, $\eps\neq 1$, then $\f_j\in {\cal C}^\eps_{\rm loc}(\Omega)$ for any $j\ge 1$.

3. If in addition $V\in{\cal C}^\infty(\Omega)$, then $\f_j\in {\cal C}^\infty(\Omega)$  for any $j\ge 1$.
\end{Lemma}

\proof 
1. The quadratic form $Q(\f)\!:=\!\displaystyle\int_{\R^n}\!V(x)|\f(x)|^2~\!dx$ satisfies
\begin{equation}
\label{eq:comp}
|Q(\f)|\le \|V\|_{\frac n{2s}}\|\f\|^2_{\sstar}
\end{equation}
by  H\"older's inequality. Hence, $Q$ is bounded in $\mathcal D^s$.  If $V\in{\cal C}^\infty_0(\R^n)$ then 
$Q$  generates a compact operator in $\mathcal D^s$ by the Rellich theorem. Since any arbitrary $V\in L^{\frac n{2s}}(\R^n)$ 
can be approximated in $L^\frac n{2s}(\R^n)$ by smooth and compactly supported functions, the corresponding operator is  compact as well, because of by (\ref{eq:comp}). 
So, the discreteness of the spectrum and the completeness of $(\f_j)$ follow by the Hilbert--Schmidt theorem. 
The equalities (\ref{eq:eigen}) hold by well known variational principle, see e.g. \cite[Sec. 10.2]{BS10}. 

Now we invoke the Green representation formula for (\ref{eq:general_eigenvalue}),
\begin{equation}
\label{eq:Green1}
\f_j(x)=C(n,s)(\Ds\f_j)*|x|^{2s-n}=C(n,s)\mu_j \irn \frac {V(\xi)\f_j(\xi)}{|x-\xi|^{n-2s}}\, d\xi~\!.
\end{equation}
Since the kernel is positive, the principal eigenfunction $\f_1$ is positive, and the corresponding eigenvalue $\mu_1$ is simple \cite{Jen}. 
On the other hand, for any $j>1$ we have
\begin{equation*}
\mu_j\irn V(x)\f_j(x)\f_1(x)~\!dx=(\Ds \f_j,\f_1)=0,
\end{equation*}
thus $\f_j$ can not have constant sign. 
 
\medskip

2. We split the integral in (\ref{eq:Green1}) into two parts:
$$
\f_j(x)=C(n,s)\mu_j \Big(
{\int\limits_{\R^n\setminus\Omega}\frac {V(\xi)\f_j(\xi)}{|x-\xi|^{n-2s}}\, d\xi+\int\limits_\Omega \frac {V(\xi)\f_j(\xi)}{|x-\xi|^{n-2s}}\, d\xi}\Big) .
$$
Since the first integral is a smooth function of $x\in\Omega$, we only have to deal with the {second} one.

We know that $\f\in L^{2^*_s}(\R^n)$. If $\eps<1$,  we use (\ref{eq:Green1}) and  statement $i)$ in Proposition \ref{Riesz} with 
$\alpha=2s$, to improve the integrability exponent for $\f_j$
which, in turns, improves the integrability exponent of $V\f_j$. A bootstrap procedure provides, in a finite number of steps, $V\f_j\in L^p(\Omega)$ for some $p>\frac n{2s}$. Then  statement
$ii)$ in Proposition \ref{Riesz} gives $\f_j\in{\cal C}^0(\Omega)$ and thus $V\f_j\in L^{\frac n{2s-\eps}}_{\rm loc}(\Omega)$. Finally, part $ii)$ in
Proposition \ref{Riesz} gives $\f_j\in{\cal C}^\eps_{\rm loc}(\Omega)$.

If $\eps>1$ then we can repeat the same steps up to obtain $\f_j\in{\cal C}^0(\Omega)$. Then we differentiate (\ref{eq:Green}), put $\alpha=2s-1$ and apply part $ii)$ in 
Proposition \ref{Riesz} to obtain $\nabla{\f}\in {\cal C}^{\eps-1}_{\rm loc}(\Omega)$, that again gives ${\f}\in {\cal C}^\eps_{\rm loc}(\Omega)$.

\medskip

3. The last claim follows from \cite[Ch. V, Theorem 4]{St} and the bootstrap argument.
\QED

\begin{Lemma}
\label{L:radial}
Let a positive weight $V\in L^{\frac n{2s-\eps}}_{\rm loc}(\R^n)$ be symmetric-decreasing.\footnote{In fact, this assumption restricts only the behavior of $V$ at zero.} For $s\le \frac 12$, assume in addition that $V\in {\cal C}^\beta_{\rm loc}(\R^n\setminus\{0\})$ with $\beta>1-2s$.
Then for any $\mu\in\R$, the problem (\ref{eq:general_eigenvalue}) has at most one linearly independent radial eigenfunction.
\end{Lemma}

\proof
We follow the outline of the proof in \cite[Theorem 1]{FLS}. 
Notice that the argument in \cite{FLS} cannot be applied directly because the weight in the right-hand side of (\ref{eq:general_eigenvalue}) 
might be singular at the origin.

We introduce the Caffarelli--Silvestre extension \cite{CS} of any function $\f\in\mathcal D^s$, that is the solution $\Phi$ of the boundary value problem
\begin{equation}
\label{eq:ACS}
-\div (y^{1-2s}\nabla \Phi)=~0\quad \mbox{in}\quad \mathbb R^n\times\mathbb R_+~,\qquad \Phi\big|_{y=0}=\f~,
\end{equation}
satisfying
$$
C_s\iirn y^{1-2s}|\nabla \Phi|^2\,dxdy= (\Ds\f,\f)
$$
for some explicitly known constant $C_s$. 
The eigenvalue problem (\ref{eq:general_eigenvalue}) can be rewritten as follows,
\begin{equation}
\label{eq:ACS1}
-C_s\cdot\lim\limits_{y\to0^+} 
y^{1-2s }\partial_y\Phi(x,y)=\mu V(x)\f(x),\qquad x\in\mathbb R^n~\!,
\end{equation}
so that
\begin{equation}
\label{eq:ARayleigh}
\mu= C_s\,\frac {\displaystyle{\iirn y^{1-2s}|\nabla \Phi|^2\,dxdy}}{\displaystyle{\irn{V(x)|\f|^2~\!dx}}}~\!.
\end{equation}
In general, problem (\ref{eq:ACS})--(\ref{eq:ACS1}) admits separation of variables; we can write its solutions in the form
\begin{equation}
\label{spherical}
\Phi(x,y)=W(r,y)Y(\Theta)~,\qquad \f(x)=h(r)Y(\Theta),
\end{equation}
where $(r,\Theta)$ are spherical coordinates in $\R^n$ and $Y$ is a spherical harmonic. 

Now we turn to the proof of the Lemma. If $\f(x)$ is a radially symmetric eigenfunction for (\ref{eq:general_eigenvalue}), then its extension 
$\Phi(x,y)$ is radially symmetric in the $x$-variable as well, and we have
$$
\Phi(x,y)=W(r,y);\qquad \f(x)=h(r)~\!.
$$

Since $V$ is positive and symmetric-decreasing, it is bounded outside of the origin. 
By \cite[Proposition B.1]{FLS}, we have $\f\in {\cal C}^{1+\delta}(\R^n\setminus\{0\})$ for some $\delta>0$. Next, Lemma \ref{L:general_spectrum} gives $\f\in{\cal C}^\eps_{\rm loc}(\R^n)$ and therefore $h\in {\cal C}^\eps_{\rm loc}(\overline{\R_+})$. So, to prove the Lemma it is sufficient to show that if $h(0)=0$ then $W\equiv0$.

We rewrite the representation formula (\ref{eq:Green1}) as follows, 
\begin{equation*}
 h(|x|)=\mu C(n,s)\irn \frac {V(|\xi|)h(|\xi|)}{|x-\xi|^{n-2s}}\, d\xi~\!.
\end{equation*}
The inclusion $h\in {\cal C}^\eps_{\rm loc}(\overline{\R_+})$ and the assumption  $h(0)=0$ reduces the order of singularity of the integrand at $\xi=0$. In turns, this gives a better H\"{o}lder estimate for $h$. Repeating this argument we obtain $h\in{\cal C}^{\beta_1}_{\rm loc}(\overline{\R_+})$ for any $\beta_1<2s$. 

Next, we rewrite  problem (\ref{eq:ACS})--(\ref{eq:ACS1}) in polar coordinates to obtain that the pair $W$, $h$ solve
$$
 \aligned
 -\, & \partial^2_{rr}W-\frac {n-1}{r}\,\partial_rW-\partial^2_{yy}W-\frac {1-2s}{y}\,\partial_yW=0, \qquad y>0;\\
 & W(r,0)=h(r);\qquad C_s\cdot\lim\limits_{y\to0^+} y^{1-2s }\partial_yW(r,y)+\mu V(r)h(r)=0
;\\
& \iirn y^{1-2s}(|\partial_r W|^2+|\partial_y W|^2)\,dxdy<\infty
~\!.  
 \endaligned
$$
Following \cite{FLS}, we introduce the function
$$
H(r)=C_s\int\limits_0^\infty y^{1-2s}\big[\big|\partial_rW|^2-|\partial_yW|^2\big]\,dy+\mu V(r)h^2(r).
$$
Proposition B.2 in \cite{FLS} gives
$$
H(\infty):=\lim_{r\to\infty} H(r)=0. 
$$
Moreover, $h(0)=0$ implies $|h(r)|\le C_\beta r^\beta$, thus $\lim\limits_{r\to0}V(r)h^2(r)=0$ because of the summability of assumption on $V$.  In addition, $\partial_rW(0,y)\equiv0$ by symmetry and thus
$$
H(0)= -\,C_s\int\limits_0^\infty y^{1-2s}\big(\partial_yW(0,y)\big)^2\,dy\le0~\!.
$$
Finally, repeating the proof of \cite[Lemma 4.1]{FLS} we conclude that 
$$
H'(r)\le -2C_s\,\frac {n-1}{r}\,\int\limits_0^\infty y^{1-2s}|\partial_rW|^2\,dy\le 0
$$
in the sense of distributions. Therefore, $H$ is non-increasing. Since $H(0)\le H(\infty)$,  we infer that  $H\equiv0$. This gives $\partial_rW\equiv0$ a.e. and hence  $W=W(t)$. But this implies $h(r)=const$, therefore $h\equiv0$ and $W\equiv0$.
 \QED

\section{Proof of Theorem \ref{T:CKN}}
\label{S:T11}

Let $u$ be a weak solution to  (\ref{eq:CKN}). As mentioned in the introduction,
 $u$ is radially symmetric about the origin and radially decreasing. 
Also, notice that $u$ solves (\ref{eq:general_eigenvalue}) for $\mu=1$ with weight $V(x)=|x|^{-bq}u^{q-2}\in L^\frac{n}{2s}(\R^n)$,
so that we can write the Green representation formula for (\ref{eq:CKN}),
\begin{equation}
\label{eq:Green}
 {u}(x)=C(n,s)\irn \frac {|\xi|^{-bq}{u}^{q-1}(\xi)}{|x-\xi|^{n-2s}}\, d\xi~\!.
\end{equation}
By repeating literally the proof of  Lemma 6 in \cite{Us} one first obtains that $u\in L^\infty(\R^n)$.
Then, using \cite[Proposition 2.6]{FW} one infers that the $s$-Kelvin transform
\begin{equation}\label{Kelvin}
x\mapsto y=\frac x{|x|^2},\qquad u(x)\mapsto \widetilde u(y)=\frac {1}{|y|^{n-2s}}\, u\Big(\frac y{|y|^2}\Big) 
\end{equation}
maps a solution of (\ref{eq:CKN}) to a solution of (\ref{eq:CKN}). This gives
\begin{equation}
\label{L:U2}
\frac {C_1}{1+|x|^{n-2s}}\le {u}(x)\le \frac {C_2}{1+|x|^{n-2s}}.
\end{equation}
Notice that the constants in (\ref{L:U2}) and in the  estimates that follow depend on the choice of ${u}$.

Thanks to Lemma \ref{L:general_spectrum}, from (\ref{L:U2}) we infer that ${u}\in {\cal C}^\infty(\R^n\setminus\{0\})\cap {\cal C}^\eps(\R^n)$ for any $\eps<2s-bq$. Thus $i)$ in Theorem \ref{T:CKN} follows by choosing $u=\bu_t$.

\medskip

Now we prove $ii)$.
Let $u$ be a solution to (\ref{eq:CKN}). Then $u$ is radially symmetric, radially decreasing and continuous on $\R^n$.
Take $t>0$ such that $u(0)=\bu_t(0)$ and put 
$$
\f(x)=u(x)-\bu_t(x)~,\qquad {\textsc V}(x)=\begin{cases}
|x|^{-bq}~\dfrac{(u(x))^{q-1}-(\bu_t(x))^{q-1}}{u(x)-\bu_t(x)}&\text{if $u(x)\neq \bu_t(x)$}\\
|x|^{-bq}~(q-1)(\bu_t(x))^{q-2}&\text{if $u(x)= \bu_t(x)$}
\end{cases}
$$
Then ${\textsc V}$ is radial and satisfies the regularity assumptions in Lemma \ref{L:radial}. It turns out that ${\textsc V}$ is symmetric-decreasing,
thanks to the next calculus lemma.

\begin{Lemma}
Let $f$ be a convex function on $\R_+$. If $u$ and $v$ are decreasing (increasing) 
positive functions on $\R_+$, then $g=\frac {f(u)-f(v)}{u-v}$ is decreasing (increasing) on $\R_+$.
\end{Lemma}
\proof
It is sufficient to assume all functions smooth. We calculate
$$
g'=\frac {u'}{(u-v)^2}\,\big(f(v)-f(u)-f'(u)(v-u)\big)+
\frac {v'}{(u-v)^2}\,\big(f(u)-f(v)-f'(v)(u-v)\big),
$$
and the statement follows. $\square$

\medskip

We can now continue the proof of the Theorem. Since ${\textsc V}$ is symmetric-decreasing and $\f$ is a radial solution of 
$$
\Ds \f={\textsc V}(x)\f~,\quad \f(0)=0,
$$
then Lemma \ref{L:radial} applies and gives $\f\equiv 0$. Thus $ii)$ is proved. 

Before proving  $iii)$ it is convenient to 
point out the next observation.

\begin{Remark}\label{R:Kelvin}
Let us notice that by $ii)$, the transform (\ref{Kelvin}) maps $\bu_t$ to $\bu_\tau$ for some $\tau>0$. 
From now on we assume that $\bu_1$ is a fixed point of the $s$-Kelvin transform.
\end{Remark}

To go further we  study in detail the action of the group of isometries $\mathcal D^s\to \mathcal D^s$ parametrized
by $t>0$ and given by
$$
%\mathcal I(t)\in \mathcal B(\mathcal D^s,\mathcal D^s)~,\qquad 
\mathcal I(t)u(x):=t^{\frac {2s-n}{2}}u\big(\frac {x}{t}\big).
$$
Notice that $\bu_t=\mathcal I(t)\bu_1$. Since $\bu_t$ is a smooth function on $\R^n\setminus\{0\}$, we can differentiate the identity $\bu_t=\mathcal I(t)\bu_1$ with respect to $t$ to  obtain
\begin{equation}
\label{eq:dot}
\budot_t=\frac 1t \mathcal I(t)\budot_1~,
\end{equation}
where the radial function $\budot_1\in {\cal C}^\infty(\R^n\setminus\{0\})$ is given by
\begin{equation}
\label{eq:budot}
\budot_1(x)=-x\cdot\nabla\bu_1(x)-\frac {n-2s}{2}\bu_1(x).
\end{equation}

For $x\ne0$, consider the integral
\begin{equation*}
\irn \underbrace{\frac {|\xi|^{-bq}\bu_1^{q-2}(\xi)|\xi\cdot\nabla\bu_1(\xi)|}{|x-\xi|^{n-2s}}}_{=:\phi(x,\xi)}\, d\xi=I_1+I_2+I_3:= 
\int\limits_{|\xi|<\frac {|x|}2}\!\!\!\!\!\!\phi ~+ \int\limits_{\frac {|x|}2<|\xi|<2|x|}\!\!\!\!\!\!\!\!\!\phi~~+\int\limits_{|\xi|> 2{|x|}}\!\!\!\!\!\!\phi.
\end{equation*}
Easily,  the integral $I_2$  converges. Furthermore, the estimate (\ref{L:U2}) implies
$$
I_1+I_3\le C(x)\int\limits_0^\infty \frac {r^{n-1}r^{1-bq}|\nabla \bu_1(r)|}{(1+r^{n-2s})^{q-1}}\,dr\stackrel{*}{\le} C(x)\int\limits_0^\infty |\nabla \bu_1(r)|\,dr~\!.
$$
(the inequality ($*$) follows from (\ref{eq:bq})).
Since $\bu_1$ is symmetric-decreasing, the last integral converges. Moreover, this convergence is uniform with respect to $x$ in any compact set bounded away from the origin.

This allows us to differentiate the equality (\ref{eq:Green}) for $u=\bu_t$ with respect to $t$. We arrive at
\begin{equation}
\label{eq:Green-dot}
 \budot_t(x)=(q-1)C(n,s)\irn \frac {|\xi|^{-bq}\bu_t^{q-2}(\xi)\budot_t(\xi)}{|x-\xi|^{n-2s}}\, d\xi~\!.
\end{equation}
We infer that $\budot_t$ is an eigenfunction to (\ref{eq:Green1}) with weight $V(x)=|x|^{-bq}\bu_t^{q-2}(x)$. By the estimate (\ref{L:U2}) we can apply part 2 of Lemma \ref{L:general_spectrum}. So, $\budot_t$ is bounded and H\"{o}lder continuous in $\R^n$. Also it is smooth outside the origin. 
Finally, the $s$-Kelvin transform gives
\begin{equation}
\label{L:U3}
|\budot_1(x)|\le \frac {C_3}{1+|x|^{n-2s}}.
\end{equation}
The estimates (\ref{L:U2}) and (\ref{L:U3}) show that $|x|^{-bq}\bu_1^{q-2}(x)\budot_1(x)\in L^{\frac{2n}{n+2s}}(\R^n)\subset (\mathcal D^s(\R^n))'$, and (\ref{eq:Green-dot}) gives $\budot_1\in \mathcal D^s(\R^n)$.

Repeating this procedure we can differentiate (\ref{eq:Green-dot}) with respect to $t$ once more. This gives the integral equation for $\buddot_t$, from which we derive, similarly to previous steps,
$$
\buddot_1\in \mathcal D^s(\R^n)\cap{\cal C}^\alpha(\R^n)\cap {\cal C}^\infty(\R^n\setminus\{0\});\qquad |\buddot_1(x)|\le \frac {C_4}{1+|x|^{n-2s}}.
$$
Since $\budot_t\in\mathcal D^s(\R^n)\setminus\{0\}$, we obtain $iii)$ in Theorem \ref{T:CKN}.
\QED

\section{\!\!\!The linearized problem and proof of Theorem \ref{T:nondegenerate}}
\label{S:T12}

Consider the functional on $\mathcal D^s$,
$$
E_0[u]=\frac 12(\Ds u,u)
-\frac 1q\int\limits_{\R^n}|x|^{-\bq q}u_+^q~\!dx,
$$
where $u_+=\max\{u,0\}$. Recalling that the truncation operator $u\mapsto u_+$ is continuous in $\mathcal D^s$ for $s\in(0,1]$, see \cite[Theorem 5.5.2/3]{RS},
and using (\ref{eq:HS_ineq}),
one can prove in a standard way that $E_0$ is of class ${\cal C}^2$, with first and second order differentials given by distributional equalities
$$
\aligned
E'_0[u] = &\ \Ds u- |x|^{-\bq q}u_+^{q-1}, \\
E''_0[u]\f = &\ \Ds\f-(q-1) |x|^{-\bq q}u_+^{q-2}\f~\!.
\endaligned
$$

\begin{Remark}
Let $u,\f,\psi\in\mathcal D^s$. The next identities for $t>0$ can be checked by elementary change of variables:
\begin{gather}
\nonumber
E_0[\mathcal I(t)u]=E_0(u);\qquad (E'_0[\mathcal I(t)u],\f)=(E'_0[u],\mathcal I(t^{-1})\f);\\
\label{eq:second}
E''_0[\mathcal I(t)u](\f,\psi)=E''_0[u](\mathcal I(t^{-1})\f,\mathcal I(t^{-1})\psi).
\end{gather}
\end{Remark}

For any $t>0$ we have that $E'_0[\bu_t]=0$ and the kernel of $E''(\bu_t)$ is  the set of
solutions to the linearized problem (\ref{eq:lin}). By the results in Section \ref{S:eigenvalue} (with weight $V(x)=|x|^{-\bq q}\bu_t^{q-2}$), the related
eigenvalue problem 
\begin{equation*}\tag{$\mathcal E_t$}
\label{eq:eigenvalue}
\Ds \f=\mu |x|^{-\bq q}\bu_t^{q-2}\f~,\qquad \f\in\mathcal D^s,
\end{equation*}
has a discrete, non decreasing sequence $(\mu_j)$ of eigenvalues that admit a variational characterization (\ref{eq:eigen}). Since the energy $E_0$ is invariant with respect to the action of the transforms $\mathcal I(t)$, 
the eigenvalues $\mu_j$ do not depend on $t>0$.

Clearly $\mu_1=1$, and 
the first eigenfunction is $\bu_t$. Next, we deal with the second eigenvalue.

\begin{Lemma}
%\label{L:second}
The eigenvalue $\mu_2$ equals $q-1$. 
\end{Lemma}

\proof
By part $iii)$ of Theorem \ref{T:CKN}, $\budot_t\in\mathcal D^s$ for any $t>0$, hence $E''_0[\bu_t]\budot_t=0$.
Thus $\mu=q-1$ is an eigenvalue for (\ref{eq:eigenvalue}), and $\mu_2\le q-1$. 

Next, recall that $\bu_t$ solves (\ref{eq:CKN}) and achieves the best constant $S_q$. Thus
$$
\irn|x|^{-\bq q}\bu_t^{q-1}\f~\!dx=(\Ds\bu_t,\f)~,\quad J''[\bu_t](\f,\f)\ge 0\quad\text{for any $\f\in\mathcal D^s$,}
$$
where
$J[u]= \dfrac {(\Ds u,u)}{\||x|^{-b} u\|_q^2}$ for $u\in\mathcal D^s\setminus\{0\}$. 
By direct computation (see for instance \cite[Lemma 3.1]{Naz2004}) we obtain
$$
\aligned
\frac {\||x|^{-\bq}\bu_t\|_q^2}{2}~J''[\bu_t](\f,\f)= &\ (\Ds\f,\f)- (q-1)~\!\irn|x|^{-\bq q}\bu_t^{q-2}\f^2~\!dx\\
+ &\ \frac {q-2}{\||x|^{-\bq}\bu_t\|_q^q}~(\Ds \bu_t,\f )^2.
\endaligned
$$
We infer that if $\f$ is orthogonal to $\bu_t$ then $(\Ds\f,\f)- (q-1)~\!\displaystyle\irn|x|^{-\bq q}\bu_t^{q-2}\f^2~\!dx\ge 0$, hence 
$\mu_2\ge q-1$. This completes the proof.
\QED

To prove Theorem \ref{T:nondegenerate} we need the following auxiliary statement.

\begin{Lemma}
%\label{L:monot}
 Let $W_0(x,y)=W_0(r,y)$ be the Caffarelli--Silvestre extension of $\bu_t$. Then $\partial_rW_0<0$ for $r>0$, $y>0$. 
 \end{Lemma}

\proof
We use the Green representation for the Caffarelli--Silvestre extension, see \cite{CS}:
$$
W_0(x,y)=C(n,s)\irn \frac {y^{2s}\bu_t(\xi)\,d\xi}{(|x-\xi|^2+y^2)^{\frac{n+2s}{2}}}.
$$
The fact that the convolution of two symmetric-decreasing functions is symmetric-decreasing is well known. We give the proof for the reader's convenience.

Since $\partial_rW_0=\sum\limits_{i=1}^nr^{-1} {x_i}~\!\partial_{x_i}\!W_0$, it suffices to prove that $\partial_{x_n}W_0(x,y)<0$ for $x_n>0$. Using the notation $x=(x',x_n)$, we derive
\begin{equation*}
\aligned
\partial_{x_n}W_0(x,y)= &\ C(n,s)(n+2s)\irn \frac {y^{2s}(\xi_n-x_n)\bu_t(\xi)\,d\xi}{(|x'-\xi'|^2+|x_n-\xi_n|^2+y^2)^{\frac{n+2s+2}{2}}}\\
= &\ C(n,s)(n+2s)\int\limits_0^\infty \int\limits_{\R^{n-1}}\frac {y^{2s}\eta\big(\bu_t(\xi',x_n+\eta)-\bu_t(\xi',x_n-\eta)\big)\,d\xi'd\eta}{(|x'-\xi'|^2+\eta^2+y^2)^{\frac{n+2s+2}{2}}}~\!.
\endaligned
\end{equation*}
From $|x_n+\eta|>|x_n-\eta|$ we infer $\bu_t(\xi',x_n+\eta)-\bu_t(\xi',x_n-\eta)<0$, since $\bu_t$ is symmetric-decreasing, and the lemma follows.
\QED

\noindent
{\bf Proof of Theorem \ref{T:nondegenerate}.}
By Lemma \ref{L:radial}, $\budot_t$ is the only radial eigenfunction corresponding to the eigenvalue $\mu_2=q-1$. Now we exclude the existence of the eigenfunctions with non-trivial spherical harmonic $Y$. In fact, we put $V(x)=|x|^{-\bq q}\bu_t^{q-2}$ in (\ref{eq:ACS})--(\ref{eq:ACS1}) and show that if $Y\not\equiv1$ in (\ref{spherical}) the quotient in (\ref{eq:ARayleigh}) is strictly greater than $q-1$. 
\medskip

Given $h(r)$, we can minimize the quotient in (\ref{eq:ARayleigh}) with respect to $Y$. This gives, modulo rotations, $Y(\Theta)=\frac{x_n}{r}$.

Since $\partial_r\bu_t<0$ and $\partial_rW_0<0$ for $r>0$, we can write
\begin{equation}
\label{eq:Apartial_r} 
 W(r,y)=\mathfrak{g}(r,y)\partial_rW_0; \qquad h(r)=g(r)\partial_r\bu_t,
\end{equation}
where $g=\mathfrak{g}|_{y=0}$.

Let $h\in{\cal C}^\infty_0(0,+\infty)$. 
Then $\mathfrak{g}$ is smooth. Using (\ref{spherical}) we rewrite (\ref{eq:Apartial_r}) as follows:
$$
\Phi(x,y)=\mathfrak{g}(r,y)\partial_{x_n}W_0;\qquad \f(x)=g(r)\partial_{x_n}\bu_t,
$$
and therefore
$$
\aligned
C_s &\iirn y^{1-2s}|\nabla \Phi|^2\,dxdy\\
=&~C_s\iirn y^{1-2s}\Big(|\nabla \partial_{x_n}W_0|^2\mathfrak{g}^2+\partial_{x_n}W_0\nabla\partial_{x_n}W_0\cdot\nabla(\mathfrak{g}^2)
+(\partial_{x_n}W_0)^2|\nabla \mathfrak{g}|^2\Big)\,dxdy.
\endaligned
$$
Integrating by parts in the second term and using the equation
$-\div (y^{1-2s}\nabla W_0)=0$ we obtain
$$
\aligned
C_s\iirn y^{1-2s}|\nabla \Phi|^2\,dxdy= &\ \irn g^2\partial_{x_n}\bu_t\Ds \partial_{x_n}\bu_t\,dx\\
+ &\ C_s\iirn y^{1-2s}(\partial_{x_n}W_0)^2|\nabla \mathfrak{g}|^2\,dxdy.
\endaligned
$$
Since $\Ds \partial_{x_n}\bu_t=\partial_{x_n}\big(|x|^{-bq}\bu_t^{q-1}\big)$, we arrive at
$$
\aligned
C_s &\iirn y^{1-2s}|\nabla \Phi|^2\,dxdy-
(q-1)\irn |x|^{-bq}\bu_t^{q-2}\f^2\,dx\\
=C_s &\iirn y^{1-2s}(\partial_{x_n}W_0)^2|\nabla \mathfrak{g}|^2\,dxdy-bq \irn g^2\bu_t^{q-1}\partial_r\bu_t\,\frac {x_n^2}{|x|^{bq+3}}\,dx.
\endaligned
$$
The right-hand side here is positive. Moreover, it is bounded away from zero if we approximate an arbitrary function $\f(x)=h(r)\frac {x_n}{r}\in \mathcal D^s(\R^n)$ by functions the supports of which are bounded and separated from the origin. This completes the proof.
\QED

\begin{Corollary}
\label{C:Morse}
There exists $\kappa>0$ independent of $t>0$ such that
\begin{equation}
\label{eq:delta}
E''_0[\bu_t](\f,\f)\ge \kappa \|\f\|_{\mathcal D^s}^2 %\irn|\Dshalf \f|^2~\!dx
\end{equation}
for any $\f\in \langle\bu_t,\budot_t\rangle^\perp$ and for any $t>0$. Moreover, the following facts hold:
\begin{itemize}
\item[$i)$] If $\f\in\mathcal D^s$ solves  $E''_0[\bu_t]\f=\gamma \Ds\budot_t$ 
for some $\gamma \in\R$, then $\f\in \langle \budot_t\rangle$, hence $\gamma =0$;
\item[$ii)$]
For any $v\in \langle \budot_t\rangle^\perp$
there exists a unique $\f\in \langle \budot_t\rangle^\perp$ such that
$E''_0[\bu_t]\f=\Ds v$. Moreover, 
$$
\kappa_*\|\f\|_{\mathcal D^s}\le \|v\|_{\mathcal D^s},
$$
where $\kappa_*=\min\{\kappa,q-2\}$.
\end{itemize}
In particular, the operator $E''_0[\bu_t]\,:\,\langle \budot_t\rangle^\perp \mapsto \Ds(\langle \budot_t\rangle^\perp)$  is isomorphism.
\end{Corollary}

\proof
We already noticed that the eigenvalues $\mu_j$ of $E''_0$ do not depend on $t>0$ because of the 
invariance of $E_0$ with respect to the transforms $\mathcal I(t)$.
Let $\mu_3$ be the third eigenvalue of (\ref{eq:eigenvalue}). Then  $\mu_3>\mu_2$ by Theorem \ref{T:nondegenerate}. Thus
(\ref{eq:delta}) holds, with $\kappa=1-\frac {q-1}{\mu_3}>0$. The last  conclusions are immediate.
\QED

\section{Proof of Theorem \ref{T:symmetry}}
\label{S:positivity}

We introduce the ${\cal C}^1$ function
$$
F(\lambda,u)=E'_0[u]+\lambda|x|^{-2s}u-(\Ds\budot_1, u)\Ds\budot_1~,\qquad F:\R\times \mathcal D^s\to (\mathcal D^s)'
$$
and notice that $F(0,\bu_1)=0$. We claim that $\partial_u F(0,\bu_1)$ is invertible. Explicitly, we have  
$$
\partial_u F(0,\bu_1)\f= E''_0[\bu_1]\f-(\Ds\budot_1, \f)\Ds\budot_1~\!, \qquad\partial_u F(0,\bu_1):\mathcal D^s\to (\mathcal D^s)'~\!.
$$
By Corollary \ref{C:Morse}, $\partial_u F(0,\bu_1)$ maps isomorphically  $\langle\budot_t\rangle^\perp$ onto $\Ds(\langle \budot_t\rangle^\perp)$. Since evidently it maps $\langle\budot_t\rangle$ onto $\Ds(\langle \budot_t\rangle)$, it isomorphically maps the space $\mathcal D^s$ onto $\Ds(\mathcal D^s)=(\mathcal D^s)'$, and the claim follows.

Thanks to the implicit function theorem, there exist $\lambda_0>0$ and a neighbourhood $\mathcal U$ of $\bu_1$
such that for any $\lambda\in(-\lambda_0,\lambda_0)$, the equation 
$F(\lambda,u)=0$ has a unique solution $u\in \mathcal U$. Of course, $u$ must be radially symmetric, precisely because of the uniqueness given by the implicit function
theorem.
To conclude the proof it suffices to show that any minimizer for $S^\lambda_q$ can be properly rescaled to obtain a function
$u_\lambda\in \mathcal U$ such that $F(\lambda,u_\lambda)=0$, provided that $\lambda>0$ is small enough 
(we already noticed that any minimizer for $S^\lambda_q$ is radially symmetric if $\lambda\le 0$).

\medskip

We start by taking any $\lambda>0$ and any minimizer $u_\lambda$ for $S^\lambda_q$. Since replacing $u_\lambda\to |u_\lambda|$ decreases the quotient in the right-hand side of (\ref{eq:lambda}), see \cite[Theorem 3]{MN-HSc}, we can assume $u_\lambda$ nonnegative. We normalize $u_\lambda$ so that it solves (\ref{eq:l_CKN}), that is 
\begin{equation}\label{eq:normalize}
\|\Dshalf u_\lambda\|_2^2+\lambda\||x|^{-s}u_\lambda\|_2^2=(S^\lambda_q)^{\frac{q}{q-2}}.
\end{equation}
Inspired by the Emden-Fowler transform,  we  introduce the functions  $v,w\,:\,\R\to\R$ given by
$$
v(\zeta)=e^{\frac{2s-n}{2}\zeta}~\bu_1(e^{-\zeta})~,\qquad w(\zeta)=
e^{\frac{2s-n}{2}\zeta}\int\limits_{\S^{n-1}}u_\lambda (e^{-\zeta}\sigma)~\!d\sigma
$$
(here we identified the radial function $\bu_1$ with a function of $r=|x|$).
 Using H\"older inequality and (\ref{eq:bq}) we obtain
$$
\int\limits_{-\infty}^\infty v^q~\!d\zeta\le c\int\limits_{-\infty}^\infty e^{(bq-n)\zeta}\int\limits_{\S^{n-1}}u_\lambda^q (e^{-\zeta}\sigma)~\!d\sigma d\zeta=
c\irn|x|^{-bq}u_\lambda^q~\!dx,
$$
(here $c$ depends only on $n, q$), that gives $w\in L^q(\R)$.

Further, by our choice in Remark \ref{R:Kelvin} we have that $v(\zeta)\equiv v(-\zeta)$; Theorem \ref{T:CKN} and formula (\ref{eq:budot})  give us  
$$
v\in {\cal C^\infty}(\R)~,\quad %v'(\zeta)=e^{\frac{2s-n}{2}\zeta}\budot_1(e^{-\zeta})~,\quad 
0\le v(\zeta)\le \frac{C_2}{\cosh (\frac{n-2s}{2}\zeta)}~\!.
$$
We infer that $\lim\limits_{|\zeta|\to\infty} v(\zeta)=0$ and $v\in L^p(\R)$ for any $p\in[1,\infty]$. In particular,  there exists $t_\lambda>0$ such that
$ t_\lambda$ achieves the maximum of the smooth function
$$
\aligned
f(t):= &\ (v^{q-1}*w)(\log t)~=\int\limits_{-\infty}^\infty (v(\zeta-\log t))^{q-1}w(\zeta)~\!d\zeta
\\
= &\ \irn t^{\frac{n-2s}{2}(q-1)}|x|^{-bq}\bu_1^{q-1}(xt)u_\lambda(x)~\!dx=\irn|x|^{-bq}\bu_{1/t}^{q-1}u_\lambda~\!dx~\!.
\endaligned
$$
Recall that $\bu_{1/t}=\mathcal I(1/t)\bu_1=(\mathcal I(t))^{-1}\bu_1$ solves (\ref{eq:CKN}) and that  
$\mathcal I(1/t)$ is isometry in $\mathcal D^s$. Thus we have
$$
\aligned
f(t)= &\ ( \Ds \bu_{{1/t}}, u_\lambda)=
(\Ds \bu_1, \mathcal I(t)u_\lambda);\\
f'(t)= &\ -\frac 1{t^2}\,(\Ds\budot_{1/t}, u_\lambda)
=-\frac 1t\,(\Ds\budot_1, \mathcal I(t)u_\lambda).
\endaligned
$$
Since $t_\lambda$ achieves the maximum of $f$, we have
$(\Ds\budot_1, \mathcal I(t_\lambda) u_\lambda )=0$ and   
$$
\aligned
(\Ds \bu_1, \mathcal I(t_\lambda)u_\lambda) =f(t_\lambda)
\ge f(\tau^{-1}t_\lambda)= &\ (\Ds \bu_{1}, \mathcal I(\tau^{-1}t_\lambda)u_\lambda)\\
= &\ \langle \Ds \bu_{\tau}, \mathcal I(t_\lambda)u_\lambda\rangle
\endaligned
$$
for any $\tau>0$.
We see that, eventually replacing $u_\lambda$ with $\mathcal I(t_\lambda)u_\lambda$, we can assume  
\begin{gather}
\label{eq:perp}
(\Ds\budot_1, u_\lambda)=0;\\
\label{eq:perp_bis}
%\|\Dshalf u_\lambda-\Dshalf \bu_1\|^2_2=\min_{\tau>0}\|\Dshalf u_\lambda-\Dshalf \bu_\tau\|^2_2~\!.
\|u_\lambda-\bu_1\|_{\mathcal D^s}=\min_{\tau>0}\|u_\lambda-\bu_\tau\|_{\mathcal D^s}~\!.
\end{gather}
From (\ref{eq:perp}) and since $u_\lambda$ solves (\ref{eq:l_CKN}) we infer that $F(\lambda,u_\lambda)=0$ for any $\lambda>0$.
To conclude the proof we only need to show that $u_\lambda\in \mathcal U$ for $\lambda$ small enough.

Take any sequence $\lambda_h\searrow 0$. 
By (\ref{eq:normalize}), $u_\lh$ is a bounded minimizing sequence for $S_q$, so we can suppose that $u_\lh$ converges weakly in $\mathcal D^s$. Arguing as in the proof of \cite[Lemma 4.2]{MN_SB}, we can rescale $u_\lh$ so that its weak limit $u$ is non-zero, hence $u$ is a (nonnegative) solution of (\ref{eq:CKN}) and $u_\lh\to u$ strongly in $\mathcal D^s$. 
Thanks to
the uniqueness result in Theorem \ref{T:CKN}, we see that there exists $\widehat\tau>0$ and a sequence $\tau_h>0$ such that 
$\mathcal I(\tau_h)u_\lh - \bu_{\widehat\tau}\to 0$ in $\mathcal D^s$. But then  (\ref{eq:perp_bis}) gives
$$
\aligned
\|u_\lh-\bu_1\|_{\mathcal D^s}\le \|u_\lh-\bu_{\widehat\tau\tau_h^{-1}}\|_{\mathcal D^s}
= &\ \|u_\lh-\big(\mathcal I(\tau_h^{-1})\bu_{\widehat\tau}\big)\|_{\mathcal D^s}\\
= &\ \|\mathcal I(\tau_h)u_\lh-\bu_{\widehat\tau}\|_{\mathcal D^s}=o(1).
\endaligned
$$
Hence $u_\lh\to \bu_1$, that is enough to conclude.
\QED

\begin{Remark}
 In fact, the assumption that $u_\lambda$ is a minimizer for (\ref{eq:lambda}) is used only to show that $u_\lambda$ are bounded for $\lambda\searrow 0$. We conjecture that there is a $\lambda_0>0$ such that for $\lambda\in (-H_s,\lambda_0]$ any nonnegative solution to (\ref{eq:l_CKN}) is radially symmetric. However, this problem is open.
\end{Remark}

\section{Dimension reduction and proof of Theorem \ref{T:pert0}}
\label{S:existence}

Given $k\in L^\infty(\R^n)$, we put
$$
G[u]=\frac 1q\irn k(x)|x|^{-bq}u_+^{q}~\!dx~\!.
$$
For any $\eps\in\R$ we introduce the energy functional on $\mathcal D^s$ given by
$$
E_\eps[u]=E_0[u]-\eps G[u]=
\frac 12\int\limits_{\R^n}|\Dshalf u|^2~\!dx-\frac 1q\int\limits_{\R^n}|x|^{-\bq q}(1+\eps k(x))u_+^q~\!dx~\!.
$$
Evidently, $E_\eps\in{\cal C}^2$, 
and any critical point $u$ for $E_\eps$ is a weak solution to 
$$
\Ds u=(1+\eps k(x))|x|^{-bq}u_+^{q-1}~\!.
$$
If $u\neq 0$ and $|\eps|\|k\|_\infty\le 1$, then $u$ is positive by the strong maximum principle  \cite{Sil}. Hence, $u$ solves (\ref{eq:eps_CKN}).

\medskip

In order to face the problem  $E'_\eps[u]=0$ for $\eps$ close to zero we combine variational methods with a Lyapunov-Schmidt technique,
in the spirit of \cite{AM}. The next lemma is the crucial step.

\begin{Lemma}[Dimension reduction]
\label{L:DR}
There exist $\eps_0> 0$ such that the problem
\begin{equation}\label{eq:perturb}
\displaystyle{E'_\eps[u]= \frac{(E'_\eps[u],\budot_t)}{\|\Dshalf \budot_t\|_2^2}\Ds\budot_t}~,\qquad u\in\langle \budot_t\rangle^\perp,
\end{equation}
has a nontrivial solution $u=U^\eps_t$ for any $(\eps,t)\in (-\eps_0,\eps_0)\times \R_+$.  Moreover, this solution is unique in a neighbourhood of $\bu_t$, the function $(\eps,t)\mapsto U^\eps_t$ is ${\cal C}^1$-smooth, and the following facts hold:
\begin{itemize}
\item[$i)$] $\|\Dshalf(U^\eps_t-\bu_t)\|_2=O(\eps)$ as $\eps\to 0$, uniformly with respect to  $t\in\R_+$;
\item[$ii)$] For any $\eps\in (-\eps_0,\eps_0)$, the curve $t\mapsto U^\eps_t$ is a natural constraint for $E_\eps$, that is, 
$$
\frac{d}{dt}E_\eps[U^\eps_{t}]\big|_{t=t_*}\!\!\!=0 \quad \text{for some}\quad t_*>0 
\qquad\Longleftrightarrow\qquad E'_\eps[U^\eps_{t_*}]=0;
$$
\item[$iii)$] Assume in addition
$\lim\limits_{x\to 0}k(x)= \lim\limits_{|x|\to \infty}k(x)=0$. Then  $U^\eps_t-\bu_t\to 0$ in $\mathcal D^s$ as $t\to 0$ and as $t\to \infty$, uniformly with respect to  $\eps$.
\end{itemize}
\end{Lemma}

\proof
We basically follow the outline of the arguments in \cite[Sections 3 and 4]{FeSc} but we considerably simplify the proofs there. Moreover, since 
the solution $\bu_1$ to the unperturbed problem is not explicitly known, in some of the steps the proof needs  more care with respect to \cite{FeSc}.   

In order to shorten formulas we denote by $\|\f\|$  the norm 
of $\f\in  \mathcal D^s$, instead of  $\|\f\|_{\mathcal D^s}$. The norm in dual space $(\mathcal D^s)'$ is denoted by $\|\cdot\|'$. Thus $\|\Ds\f\|'=\|\f\|$ for any $\f\in\mathcal D^s$.

If $X, Y$ are Banach spaces, we denote by $|\!|\!|\cdot |\!|\!|_{X\to Y}$ the standard norm in $\mathcal B(X,Y)$, which is the space of linear and continuous operators $X\to Y$. If $X$ and $Y$ are clear by definition, we write simply $|\!|\!|\cdot |\!|\!|$. For instance, if $J:\mathcal D^s\to\R$ is a smooth functional, then for any $u\in \mathcal D^s$ we have
$J'[u]\in (\mathcal D^s)'$, $J''[u]\in \mathcal B(\mathcal D^s,(\mathcal D^s)')$ and we have
$$
\|J'[u]\|'=\sup_{\f\in\mathcal D^s\atop \|\f\|=1}|(J'[u],\f)|~,\qquad 
|\!|\!|J''[u]|\!|\!|=\sup_{\psi\in\mathcal D^s\atop \|\psi\|=1}\|J''[u]\psi\|'=\sup_{\psi\in\mathcal D^s\atop \|\psi\|=\|\f\|=1}|J''[u](\psi,\f)|.
$$

We introduce also the extended space 
$$
{\mathcal D}^s_{\!\times}=\mathcal D^s\times\R\quad\text{with  norm}\quad
\|(\eta,\gamma)\|_{\!\times}^2:=\|\eta\|^2+\gamma^2~\!,
$$
and its dual space $({\mathcal D}^s_{\!\times})'=(\mathcal D^s)'\times\R$, with  norm $\|(\Ds \eta,\gamma)\|'_{\!\times}=\|(\eta,\gamma)\|_{\!\times}$.

\medskip

Consider the map $\mathfrak{F}=[\mathfrak{F}_1,\mathfrak{F}_2]^\top:(\R\times\R_+)\times{\mathcal D}^s_{\!\times}\to ({\mathcal D}^s_{\!\times})'$:
$$
\aligned
\mathfrak{F}_1(\eps,t;\eta,\gamma):= &\ E'_\eps[\bu_t+\eta]+t\gamma\Ds\budot_t\in(\mathcal D^s)'~,\\
\mathfrak{F}_2(\eps,t;\eta,\gamma):= &\ t(\Ds\budot_t,\eta)\in\R
\endaligned
$$
(the multiplier $t$ in both entries is a normalization factor; notice  that 
$t  \| \budot_t\| = \|   \budot_1 \|$ does not depend on $t$ by (\ref{eq:dot})).

The function $\mathfrak{F}$ is continuously differentiable (for the derivative with respect to $t$ use part $iii)$ in Theorem \ref{T:CKN}) and  $\mathfrak{F}(0,t;0,0)\equiv0$. We fix $t>0$ and solve the equation $\mathfrak{F}(\eps, t;\eta,\gamma)=0$
in a neighbourhood of $(0,t;0,0)$. To this goal we define  
$$
\mathfrak{L}(\eps,t;\eta,\gamma):=\partial_{(\eta,\gamma)}\mathfrak{F}(\eps,t;\eta,\gamma)\in {\mathcal B({\mathcal D}^s_{\!\times}, ({\mathcal D}^s_{\!\times})')};\qquad 
\mathfrak{L}(t):=\mathfrak{L}(0,t;0,0).
$$
In  matrix form, we have
$$
\mathfrak{L}(\eps,t;\eta,\gamma)=
\begin{bmatrix}
E''_\eps[\bu_t+\eta] & t\Ds\budot_t \\
t\Ds\budot_t & 0
\end{bmatrix};
\qquad
\mathfrak{L}(t)=
\begin{bmatrix}
E''_0[\bu_t] & t\Ds\budot_t \\
t\Ds\budot_t & 0
\end{bmatrix}.
$$

First, we claim that $\mathfrak{L}(t)$ is invertible, and the norm of $\mathfrak{L}(t)^{-1}$ admits the estimate independent of $t$. Indeed, for  $\f\in \mathcal D^s$ and $\zeta\in\R$ we have
$$
\mathfrak{L}(t)\!\!
\begin{bmatrix}
\f \\
\zeta
\end{bmatrix} =
\begin{bmatrix}
E''_0[\bu_t]\f+t\zeta\Ds\budot_t \\
t(\Ds\budot_t,\f)
\end{bmatrix} .
$$
Assume $\mathfrak{L}(t)\,\![\f,\zeta]^\top=0$. The vanishing of the first entry implies that $\f\in \langle \budot_t\rangle$ and $\zeta=0$ by part $i)$ in Corollary \ref{C:Morse}; on the other hand, the vanishing of the second entry gives $\f\in \langle\budot_t\rangle^\perp$. So, $\mathfrak{L}(t)$ is injective.

To prove that $\mathfrak{L}(t)$ is surjective we take $v\in\mathcal D^s$, $\gamma\in\R$, and seek for $\f\in\mathcal D^s$, $\zeta\in\R$ such that
\begin{equation}
\label{eq:surjective1}
E''_0[\bu_t]\f=-t\zeta\Ds +\Ds v~,\qquad
%\label{eq:surjective2}
t(\Ds\budot_t,\f)=\gamma.
\end{equation}
We choose
$$
\zeta= \frac{(\Ds \budot_t,v)}{t\|\budot_t\|^2},
$$
so that $v-t\zeta\budot_t\in \langle \budot_t\rangle^\perp$. By part $ii)$ in Corollary \ref{C:Morse}, we find a unique $\f^\perp \in \langle \budot_t\rangle^\perp$ such that
$$
E''_0[\bu_t]\f^\perp=\Ds(v-t\zeta \budot_t),
$$
It is easy to check that the (unique) solution to (\ref{eq:surjective1}) is given by  
$$
\f= \frac {\gamma}{t\|   \budot_t \|^2 }~\budot_t+\f^\perp~\!.
$$
We recall that $t  \| \budot_t\| = \|   \budot_1 \|$ and use part $ii)$ in Corollary \ref{C:Morse} to infer
\begin{equation*}
%\label{eq:future}
|\zeta|\le \frac{\|v\|}{\|\budot_1\|} ,\qquad \|\f\|^2\le \frac{\|v\|^{2}}{\kappa_*^2}+\frac{\gamma^2}{\|\budot_1\|^2}~\!.
\end{equation*}
Thus,
\begin{equation}\label{eq:inverse}
|\!|\!|\mathfrak{L}(t)^{-1}|\!|\!|\le c_*:=\frac 1{\min\{\kappa_*,\|\budot_1\|\}}, 
\end{equation}
and the claim follows.

Thanks to the implicit function theorem, for any $t>0$ and any $\eps$ close to zero the equation $\mathfrak{F}(\eps, t;\eta,\gamma)=0$ is uniquely solvable in a neighbourhood of $(0,t;0,0)$. We denote this solution by $[\eta^\eps_t,\gamma^\eps_t]^\top$ and put
$$
U^\eps_t:= \bu_t+\eta^\eps_t~\!.
$$
The equality $\mathfrak{F}_2(\eps, t;\eta^\eps_t,\gamma^\eps_t)=0$ gives $\eta^\eps_t\in \langle\budot_t\rangle^\perp$ and thus $U^\eps_t\in \langle\budot_t\rangle^\perp$. Further, the equality $\mathfrak{F}_1(\eps, t;\eta^\eps_t,\gamma^\eps_t)=0$ reads
$$
E'_\eps[U^\eps_t]=-t\gamma^\eps_t\Ds\budot_t.
$$  
Testing this equation with $\budot_t$ we see that it solves (\ref{eq:perturb}). The ${\cal C}^1$ regularity of the function $(\eps,t)\mapsto (U^\eps_t,\gamma^\eps_t)$ is given by the implicit function theorem.
\medskip

To prove $i)$ we need some estimates. We begin with 
\begin{equation}\label{eq:est}
|\!|\!|\mathfrak{L}(\eps,t;\eta,\gamma)-\mathfrak{L}(t)|\!|\!|%_{{\mathcal D}^s_{\!\times}\to ({\mathcal D}^s_{\!\times})'}
\le
|\!|\!|E''_0[\bu_t+\eta]-E''_0[\bu_t]|\!|\!|%_{{\mathcal D}^s\to ({\mathcal D}^s)'}
+|\eps|~\!|\!|\!|G''[\bu_t+\eta]|\!|\!|%_{{\mathcal D}^s\to ({\mathcal D}^s)'}
. 
\end{equation}
We define
$$
C_0(\rho):= \sup_{\| \eta \|\le\rho} |\!|\!|E''_0[\bu_t+\eta]-E''_0[\bu_t]|\!|\!| 
$$
and notice that $C_0(\rho)\to 0$ as $\rho\to 0$, because $E_0$ is of class ${\cal C}^2$. Moreover, $C_0(\rho)$ does not depend on $t$. Indeed, since $\mathcal I(t)$ is an isometry in $\mathcal D^s$, the relation (\ref{eq:second}) gives
$$
C_0(\rho)=\sup_{|\!|\!|\mathcal I(t^{-1})\eta|\!|\!|\le\rho} |\!|\!|E''_0[\mathcal I(t)(\bu_1+\mathcal I(t^{-1})\eta)]-E''_0[\mathcal I(t)\bu_1]|\!|\!|=
\sup_{\| \eta \|\le\rho} |\!|\!|E''_0[\bu_1+\eta]-E''_0[\bu_1]|\!|\!|.
$$
Thus we can  fix a small $\rho_0>0$ such that if $\|\eta \|\le\rho_0$ then the first term in the right-hand side of (\ref{eq:est}) does not exceed $\frac 1{3c_*}$, where $c_*$ is defined in (\ref{eq:inverse}).

Further, by the H\"older inequality and (\ref{eq:HS_ineq}) we obtain for $\|\eta \|\le\rho_0$ 
$$
|\!|\!|G''[\bu_t+\eta]|\!|\!|
\le (q\!-\!1)\|k\|_\infty\!\!\sup_{\f,\psi\in\mathcal D^s\atop\|   \f \|, \|\psi \|=1}~\irn
\!\!\!|x|^{-bq}|\bu_t+\eta|^{q-2}|\f||\psi|~\!dx\le c_1\| \bu_t+\eta \|^{q-2}\le c_2 
$$
where $c_2$ does not depend on $t$. Therefore, there is $\eps_0$ independent of $t$ such that for $|\eps|<\eps_0$ and $\|\eta \|\le\rho_0$ the second term in the right-hand side of (\ref{eq:est}) also does not exceed $\frac 1{3c_*}$. 

By the Banach inverse mapping theorem, for any $t>0$, $\gamma\in\R$, $|\eps|<\eps_0$ and $\|\eta \|\le\rho_0$ the operator $\mathfrak{L}(\eps,t;\eta,\gamma)$ is invertible, and 
$$
|\!|\!|\mathfrak{L}(\eps,t;\eta,\gamma)^{-1}|\!|\!|
=|\!|\!|\mathfrak{L}(t)^{-1}\big(\mathbb I+(\mathfrak{L}(\eps,t;\eta,\gamma)-\mathfrak{L}(t))\mathfrak{L}(t)^{-1}\big)^{-1}|\!|\!|
\le 3c_*.
$$
We are allowed to differentiate the implicit function and obtain
\begin{equation}\label{eq:d-eps}
\partial_\eps 
\begin{bmatrix}
\eta^\eps_t \\
\gamma^\eps_t
\end{bmatrix} =
-\mathfrak{L}(\eps,t;\eta^\eps_t,\gamma^\eps_t)^{-1}\partial_\eps \mathfrak{F}(\eps,t;\eta^\eps_t,\gamma^\eps_t)=\mathfrak{L}(\eps,t;\eta^\eps_t,\gamma^\eps_t)^{-1}
\begin{bmatrix}
G'[\bu_t+\eta^\eps_t] \\
0
\end{bmatrix}.
\end{equation}
Using again the H\"older inequality and (\ref{eq:HS_ineq}) we get for $\|\eta\|\le\rho_0$
\begin{equation}\label{eq:G'}
\|G'[\bu_t+\eta]\|'
\le\|k\|_\infty\sup_{\f\in\mathcal D^s\atop\|   \f \|=1}~\irn
|x|^{-bq}|\bu_t+\eta|^{q-1}|\f|~\!dx\le c_3\| \bu_t+\eta \|^{q-1}\le c_4 
\end{equation}
with $c_4$ independent of $t$. Therefore, the relation (\ref{eq:d-eps}) gives $\|\partial_\eps \eta^\eps_t\|\le c_5:=3c_*c_4$
which implies $\|U^\eps_t-\bu_t\|=\|\eta^\eps_t\|\le c_5\eps$. Thus, $i)$ is proved. 

Reducing $\eps_0$ if needed we arrive at $c_5\eps_0\le\rho_0$. Now $\eta^\eps_t$ (and thus $U^\eps_t$) is well-defined in the whole strip $t>0$, $|\eps|<\eps_0$.
\medskip

To prove $ii)$ we test (\ref{eq:perturb}) with $\dot{U}^\eps_t$. This gives
$$
\aligned
\frac{d}{dt}E_\eps[U^\eps_t]=
(E'_\eps[U^\eps_t],\dot{U}^\eps_t)= &\ 
\frac{(E'_\eps[U^\eps_t],\budot_t)}{\|\budot_t\|^2}~(\Ds\budot_t,\dot{U}^\eps_t)\\
= &\ 
\Big(1+ \frac{(\Ds\budot_t,\dot{U}^\eps_t-\budot_t)}{\|\budot_t\|^2}\Big)~\!(E'_\eps[U^\eps_t],\budot_t).
\endaligned
$$
Part $iii)$ in Theorem \ref{T:CKN} allows us to write
$$
( \Ds\budot_t,\dot{U}^\eps_t)=\frac{d}{dt}( \Ds\budot_t,U^\eps_t)-( \Ds\buddot_t,U^\eps_t)
=-( \Ds\buddot_t,U^\eps_t)
$$
(for the last equality use $U^\eps_t\in\langle \budot_t\rangle^\perp$) and, in a similar way,
 $$
( \Ds\budot_t,\budot_t)= \frac{d}{dt} ( \Ds\budot_t,\bu_t)-( \Ds\budot_t,\budot_t)  =-   ( \Ds\buddot_t,\bu_t)~\!.
 $$
We  differentiate with respect to $t$  the identity $t\budot_t = \mathcal I(t)\budot_1$ to get 
$t^2\buddot_t=\mathcal I(t)(\buddot_1-\budot_1)$. Since $\mathcal I$ is an isometry, we infer 
$$
\frac{ \big|(\Ds\budot_t,\dot{U}^\eps_t-\budot_t)\big|}{\|\budot_t\|^2}=
\frac{ \big|(\Ds\buddot_t,U^\eps_t-\bu_t)\big|}{\|\budot_t\|^2}\le \frac{t^{-2}\|\buddot_1-\budot_1\|}{t^{-2}\|\budot_1\|^2}~\!\|U^\eps_t-\bu_t\|
\le c_6\eps 
$$
by $i)$, with $c_6$ independent of $t>0$. 
Therefore, if $c_6\eps_0<1$ then 
$$
\frac{d}{dt}E_\eps[U^\eps_t]= 0\qquad \Longleftrightarrow \qquad
(E'_\eps[U^\eps_t],\budot_t)=0,
$$ 
and the latter relation is equivalent to $E'_\eps[U^\eps_t]=0$ by (\ref{eq:perturb}).
Thus  $ii)$ in the statement holds.

\medskip

To prove $iii)$ we sharpen the estimate (\ref{eq:G'}). The assumtpions on $k$ imply that the function  
\begin{equation}
\label{eq:g}
g(t):= \Big(\irn |k(tx)| |x|^{-bq}\bu_1^q~\!dx\Big)^{\frac {q-1}{q}}= \Big(\irn |k(x)| |x|^{-bq}\bu_t^q~\!dx\Big)^{\frac {q-1}{q}}
\end{equation}
is bounded, continuous  (use part $iii)$ in Theorem \ref{T:CKN}), and satisfies 
$$
\lim\limits_{t\to 0}g(t)= \lim\limits_{t\to \infty}g(t)= 0
$$
by the Lebesgue dominated convergence theorem. 

For $\|\eta\|\le\rho\le\rho_0$ we write
$$
\aligned
\vphantom{a} &
\|G'[\bu_t+\eta]\|'\le\sup_{\f\in\mathcal D^s\atop\| \f \|=1}~\irn
|k(x)||x|^{-bq}|\bu_t+\eta|^{q-1}|\f|\,dx\le I_1+I_2\\
:= &\ \sup_{\f\in\mathcal D^s\atop\| \f \|=1}~\irn
|k(x)||x|^{-bq}|\bu_t|^{q-1}|\f|\,dx+c\sup_{\f\in\mathcal D^s\atop\| \f \|=1}~\irn
|k(x)||x|^{-bq}|\bu_t+\eta|^{q-2}|\eta||\f|\,dx.
\endaligned
$$

We change the variable, use H\"older inequality and (\ref{eq:HS_ineq}) once again and arrive at
$$
I_1=\sup_{\f\in\mathcal D^s\atop\| \f \|=1}~\irn|k(tx)|~\!|x|^{-bq}\bu_1^{q-1}|\mathcal I(t^{-1})\f|~\!dx\le c_7g(t).
$$ 
In a similar way we get
$$
I_2\le c\|z_t+\eta\|^{q-2}\|\eta\|\le c_8\|\eta\|
$$
(here $c_7$ and $c_8$ do not depend on $t$). Therefore, (\ref{eq:d-eps}) gives 
$$
\|\partial_\eps \eta^\eps_t\|\le c_9\big(g(t)+\rho\big),\qquad c_9:=3c_*\max\{c_7,c_8\},
$$
which implies $\|\eta^\eps_t\|\le c_9\big(g(t)+\rho\big)\eps$. Thus, we obtain the implication
\begin{equation}\label{eq:g(t)}
\|\eta^\eps_t\|\le \rho\le\rho_0 \qquad\Longrightarrow\qquad \|\eta^\eps_t\|\le c_9\big(g(t)+\rho\big)\eps_0. 
\end{equation}
Reducing $\eps_0$ if needed we arrive at $c_9\eps_0<1$. Then (\ref{eq:g(t)}) yields
$$
\|U^\eps_t-\bu_t\|=\|\eta^\eps_t\|\le \frac {c_9\eps_0}{1-c_9\eps_0}\,g(t),
$$
and $iii)$ follows. The proof is complete.
\QED

\noindent
{\bf Proof of Theorem \ref{T:pert0}.} As in the previous proof,  we let $\|\cdot\|$  be the norm in 
$\mathcal D^s$.

Up to multiplication of $u$ by a proper constant we can assume without restriction that 
$\lim\limits_{x\to 0}k(x)= \lim\limits_{|x|\to \infty}k(x)=0$.  

\medskip

Let $U^\eps_t$ be the function given by Lemma \ref{L:DR} and write
$$
E_\eps[U^\eps_t]= E_0(\bu_t)+
\frac 12\,(\|U^\eps_t\|^2-\|\bu_t\|^2)
-\frac 1q\irn|x|^{-\bq q}(1+\eps k(x))((U^\eps_t)_+^q-\bu_t^q)\,dx
- \eps G(\bu_t)~\!.
$$
Recall that $\|\bu_t\|=\|\bu_1\|$ does not depend on $t$. 
From the statement $iii)$ in Lemma \ref{L:DR} we infer that $U^\eps_t$ is    
uniformly bounded in $\mathcal D^s$, $\|U^\eps_t-\bu_t\|=o(1)$ as $t\to 0$ and as $t\to\infty$ and therefore
$$
\big|\|U^\eps_t\|^2-\|\bu_t\|^2\big|\le (\|U^\eps_t\|+\|\bu_t\|)~\!\|U^\eps_t-\bu_t\|=o(1).
$$ 
Moreover, 
$\|~\!|x|^{-b}((U^\eps_t)_+-\bu_t)\|_q\le \|~\!|x|^{-b}(U^\eps_t-\bu_t)\|_q=o(1)$
by (\ref{eq:HS_ineq}).
Using also H\"older inequality we plainly infer
$$
\big|\irn|x|^{-\bq q}(1+\eps k(x))((U^\eps_t)_+^q-\bu_t^q)\,dx\big|\le 
c \irn|x|^{-\bq q}((U^\eps_t)_+^{q-1}+\bu_t^{q-1})\,|(U^\eps_t)_+-\bu_t|\,dx=o(1).
$$
Finally, we already noticed that
$|G(\bu_t)|\le g(t)^{q'}=o(1)$, where $g$ is the function in (\ref{eq:g}), and we can conclude that
$$
\phi^\eps(t):= E_\eps[U^\eps_t]= E_0(\bu_t)+o(1)=E_0(\bu_1)+o(1)\qquad \mbox{as}\ \  t\to 0\ \ \mbox{and as}\ \  t\to \infty.
$$
Thus,  $\phi^\eps$ has at least one critical point $t_\eps$ (in fact, $\phi^\eps$ might be constant). Hence $ U^\eps_{t_\eps}$ is a critical point for $E_\eps$ by the statement $ii)$ in Lemma \ref{L:DR}. The conclusion follows.
\QED

\begin{Remark}
Theorem 1.4 in \cite{FeSc} can be  extended to the fractional case as well, with minor modifications in the proof. Moreover, the above arguments apply to more general problems of the form
$$
\Ds u+\lambda|x|^{-2s} u
=|x|^{-bq}u^{q-1}+\eps f(x,u)~\!,
$$
where the perturbation term $f:\R^n\times\R\to \R$ is a Carath\'eodory function satisfying suitable regularity and growth assumptions at $0$ and at $\infty$.
\end{Remark}

\small
%\section{references}
%\label{References}

\end{document}